\documentclass[11pt]{article}
\usepackage{amsmath}
\usepackage{amsfonts}
\newtheorem{theorem}{Theorem}
\newtheorem{lemma}[theorem]{Lemma}
\newtheorem{proposition}[theorem]{Proposition}

\newtheorem{remark}[theorem]{Remark}
\newenvironment{proof}[1][Proof]{
\par\noindent{\em #1}. }{\hfill\framebox(6,6)\par\medskip}

\newcommand{\abs}[1]{\lvert#1\rvert}
\newcommand{\Irr}{\mathop{\rm Irr}\nolimits}
\newcommand{\Fr}{\mathop{\rm Fr}\nolimits}
\newcommand{\ep}{\varepsilon}
\renewcommand{\phi}{\varphi}

\newcommand{\om}{\omega}

\begin{document}

\title{On an
asymptotic behavior of elements
of order $p$\\ in irreducible representations
of the classical algebraic groups\\
with large enough highest weights}

\author{I. D. Suprunenko\\
Institute of Mathematics,  National Academy of Sciences of Belarus\\
  Surganov str. 11, Minsk, 220072, Belarus\\
suprunenko@im.bas-net.by}

\date{}
\maketitle

\begin{abstract}
The behavior of the images of a fixed element of order $p$ in
irreducible representations of a classical algebraic group in
characteristic $p$ with highest weights large enough with respect
to $p$ and this element is investigated. More precisely, let $G$
be a classical algebraic group of rank $r$ over an algebraically
closed field $K$ of characteristic $p>2$. Assume that an element
$x\in G$ of order $p$ is conjugate to that of an algebraic group
of the same type and rank $m<r$ naturally embedded into $G$. Next,
an integer function $\sigma_x$ on the set of dominant weights of
$G$ and a constant $c_x$ that depend only upon   $x$, and a
polynomial $d$ of degree one are  defined.  It is proved that the
image of $x$
 in the irreducible representation of $G$ with highest weight $\omega$
contains more than $d(r-m)$ Jordan blocks of size $p$ if $m$ and
$r-m$ are not too small and $\sigma_x(\omega)\geq p-1+c_x$.
\end{abstract}

\maketitle

Asymptotic lower estimates for the number of Jordan blocks of size
$p$ in the images of a fixed element of order $p$ in irreducible
representations of a classical algebraic group in characteristic
$p$ with highest weights large enough with respect to $p$ and this
element are obtained. More precisely, let $G$ be a classical
algebraic group of rank $r$ over an algebraically closed field $K$
of characteristic $p>2$. Assume that an element $x\in G$ of order
$p$ is conjugate to that of an algebraic group $G_m$ of the same
type and rank $m<r$ naturally embedded into $G$. Set
\begin{equation*}
d(r-m)=
\begin{cases}
r-m & \text{for $G=A_r(K)$},\\
2r-2m & \text{for other types.}
\end{cases}
\end{equation*}
Let $\Delta_x$ be the labelled Dynkin diagram of the conjugacy
class containing $x$ in the sense of Bala and Carter \cite{BC} and
let  $c_x$ be the sum of the labels at $\Delta_x$ for $G\neq
A_r(K)$ and the half of this sum for $G=A_r(K)$. For brevity,
throughout the article we refer to $\Delta_x$  as the labelled
Dynkin diagram of $x$. Next, an integer function $\sigma_x$ on the
set of dominant weights of $G$ that depends only upon $\Delta_x$
is defined. For $p$-restricted weights $\sigma_x$ coincides with
the canonical homomorphism determined by $\Delta_x$.  It is proved
that the image of $x$ in the irreducible representation of $G$
with highest weight $\om$ contains more than $d(r-m)$ Jordan
blocks of size $p$ if $m$ and $r-m$ are not too small and
$\sigma_x(\om)\geq p-1+c_x$.

We need some more notation to formulate the main results. Let $\om_i$
and $\alpha_i$ be  the fundamental weights and the simple roots of
$G$ (with respect to a fixed maximal torus $T$) labelled as in
\cite{B1}. Denote by $\delta_i$, $1\leq i\leq r$, the label on
$\Delta_x$ corresponding to its $i$th node. We have $0\leq
\delta_i\leq2$.  In what follows $\mathbb{Z}$ is the set of integers,
$\mathbf{X}=\mathbf{X}(G)$ is the set of weights of $G$,
$\mathbf{X^+}\subset\mathbf{X}$ is the set of dominant weights,
$\Irr=\Irr(G)$ is the set of irreducible rational representations
of $G$ (considered up to the equivalence) and $\om(\phi)$ is the
highest weight of a representation $\phi$. There exists a uniquely
determined homomorphism $\tau_x: \mathbf{X}\rightarrow\mathbb{Z}$
such that $\tau_x(\alpha_i)=\delta_i$. The weight $\mu\in\mathbf{X^+}$
is called $p$-restricted if $\mu=\sum^r_{i=1}a_i\om_i$ with all
$a_i<p$. Each weight $\om\in\mathbf{X^+}$ can be represented in
the form $\sum^s_{j=0}p^j\om^j$ where $\om^j$ are $p$-restricted.
Set $\sigma_x(\om)=\tau_x(\sum^s_{j=0}\om^j)$. Now we can state our
main result.
\begin{theorem}\label{mt} Let $\phi\in\Irr$ and $\sigma_x(\om(\phi))\geq p-1+c_x$.
 Assume that $m>1$ for $G=B_r(K)$ or $D_r(K)$, $r-m>1$ for $G=A_r(K)$,
and $>3$ for $G=B_r(K)$ or $D_r(K)$. Then the element $\phi(x)$ has
  more than $d(r-m)$ Jordan blocks of size $p$. \end{theorem}
Proposition~\ref{mp} below shows that one cannot weaken the
inequality for $\sigma_x(\om(\phi))$ in Theorem~\ref{mt} and that
the estimates obtained are asymptotically exact.
\begin{proposition}\label{mp} Let $\phi\in\Irr$ and $\om=\om(\phi)=a\om_1$ with $a<p$.
Assume that $m$ and $r$ are such as in Theorem~\ref{mt} and $x$ is
a regular unipotent element in $G_m$. Set $c=m$ for $G=A_r(K)$,
$2m$ for $G=B_r(K)$, $2m-1$ for $G=C_r(K)$, and $2m-2$ for
$G=D_r(K)$. Suppose that $p>c$. Then $\abs{x}=p$,
$\sigma_x(\om)=ac$ and $c_x=c$.   There exist constants
$N_G(a,m,p)$ and $Q_G(a,m,p)$ that depend upon the type of $G$,
$a$, $m$, and $p$ and do not depend upon $r$ such that $\phi(x)$
contains at most $N_G(a,m,p)$ Jordan blocks of size $p$ if
$p<ac<p+c-1$ and at most $d(r-m)+Q_G(a,m,p)$ such blocks if
$ac=p+c-1$. \end{proposition}

 Put $l=[(m+2)/2]$ for $G=A_r(K)$ and $l=m$
otherwise. By Lemma~\ref{lc} below,
$c_x=\sum^l_{i=1}\delta_i$ and $c_x\leq p-1$. Hence
$c_x\leq2l$.

For $\varphi\in\Irr$ define the weight $\bar{\omega}(\varphi)$ as
follows: write down the $p$-adic expansion for the weight
$\omega=\omega(\varphi)$ considered before the statement of
Theorem~\ref{mt} and set
$\bar{\omega}(\varphi)=\sum^s_{j=0}\om^j$. So
$\sigma_x(\omega(\varphi))=\tau_x(\bar{\omega}(\varphi))$.

 The study of an asymptotic behavior of elements of order $p$ in
representations of the classical groups in characteristic $p$ was
begun by the author in \cite{Supl} where a notion of a $p$-large
representation was introduced. In our present notation a
representation $\phi\in\Irr$ is $p$-large if and only if
$\sigma_x(\om(\phi))\geq p$ for a long root element $x\in G$ (an
equivalent definition from \cite{Supl}: the value of ${\bar
\om}(\phi)$ on the maximal root is $\geq p$). The common goal of
 \cite{Supl} and the present article is
 to investigate the behavior of elements of order $p$ in irreducible
representations in characteristic $p$ for a fixed $p$ and
$r\rightarrow\infty$ and to discover asymptotic regularities which
are specific for prime characteristics but do not (or almost do
not) depend upon $p$. Such properties can find applications in
recognizing representations and linear groups. According to
~\cite[Theorem~1.1]{Supl}, the image of any element of order $p$
in a $p$-large representation has at least $f(r)$ Jordan blocks of
size $p$ where
\begin{equation*}
f(r)=
\begin{cases}
2r-2 & \text{for $G=A_r(K)$},\\ 6r-7 & \text{for $G=B_r(K)$,
$p=3$},\\ 8r-10 & \text{for $G=B_r(K)$, $p>3$},\\ 4r-4 & \text{for
$G=C_r(K)$},\\ 4r-8 & \text{for $G=D_r(K)$, $p=2$},\\ 6r-10 &
\text{for $G=D_r(K)$, $p=3$},\\ 8r-16 & \text{for $G=D_r(K)$,
$p>3$}.
\end{cases}
\end{equation*}
In \cite[Theorem~1.3]{Supl} for $p>2$ and all  types of the
classical groups examples of representations $\phi\in\Irr$ such
that $\phi(x)$ has only one Jordan block of size $p$ for a long
root element $x$ and $\sigma_x(\om(\phi))=p-1$ are given. It is
also shown \cite[Theorem~1.4] {Supl} that the estimates in
\cite[Theorem~1.1]{Supl} are asymptotically exact for the groups
of type $A$, $B$, and $D$ provided $p>2$ in the last two cases.

In what follows $\mathbb{C}$ is the complex field,
$G_{\mathbb{C}}$ is  the simple algebraic group over  $\mathbb{C}$
of the same type and rank as $G$, and $\Irr_{\mathbb{C}}$ is the set
of irreducible rational representations of $G_{\mathbb{C}}$
(considered up to the equivalence).  For $\rho\in\Irr$ or
$\Irr_{\mathbb{C}}$  and a unipotent element $z\in G$ or
$G_{\mathbb{C}}$ denote by  $k_{\rho}(z)$ the degree of the
minimal polynomial of $\rho(z)$. It is well known that $k_{\rho}(z)$ is equal
to the maximal size of a Jordan block of $\rho(z)$. If $\phi\in\Irr$, then
$\phi_{\mathbb{C}}$ is   the irreducible representation of
$G_{\mathbb{C}}$ with highest weight ${\bar \om}(\phi)$. For
unipotent $x\in G$ put
$k_{\phi_{\mathbb{C}}}(x)=k_{\phi_{\mathbb{C}}}(y)$ where $y\in
G_{\mathbb{C}}$ is a unipotent element with the labelled Dynkin
diagram   $\Delta_x$ (this is correctly determined). Now let
$\abs{x}=p$. By the results \cite[Theorem 1.1,
Lemma~2.5, and Proposition~2.12]{Sump},
$k_{\phi_{\mathbb{C}}}(x)=\sigma_x(\om(\phi))+1$ and
$k_{\phi}(x)=\min\{p, k_{\phi_{\mathbb{C}}}(x)\}$. Hence if $z$ is
a long root element,  then $k_{\phi}(z)=k_{\phi_{\mathbb{C}}}(z)$
if and only if $\phi$ is not $p$-large. The results of \cite{Sump}
imply that for not very small $p$ and $r$ there exists a wide
class of representations $\phi\in\Irr$ such that
$k_{\phi}(z)=k_{\phi_{\mathbb{C}}}(z)<p$ for a long root element
$z\in G$, but $k_{\phi}(x)=p<k_{\phi_{\mathbb{C}}}(x)$ for many
other elements $x\in G$ of order $p$. In this connection in
\cite[Section 2]{SuKo} a notion of a $p$-large representation for
a given
 element $x$ of order $p$ was introduced. A representation $\phi\in\Irr$ was
called $p$-large for $x$ if $\sigma_x(\om(\phi))\geq p$. It has
been conjectured (\cite[Conjecture 1]{SuKo} that if $x\in G_m$,
$r$ is large enough with respect to $m$ and $\phi$ is $p$-large
for $x$, then $\phi(x)$ has at least $F(r)$ blocks of size $p$
where $F$ is an increasing function. Our Proposition~\ref{mp}
formally disproves this conjecture, but Theorem~\ref{mt} actually
proves a refined version of it with a stronger assumption on
$\sigma_x(\om(\phi))$. Thus for arbitrary elements $x$ there is a
gap between the class of representations $\phi\in\Irr$ with
$k_{\phi}(x)=k_{\phi_{\mathbb{C}}}(x)$ and that of representations
where asymptotic estimates for the number of Jordan blocks of size
$p$ in $\phi(x)$ hold. Perhaps for some classes of elements of
order $p$ stronger estimates than those of Theorem~\ref{mt} are
possible, but now it is not clear how to determine such classes.

The case $p=2$ is not considered here, but in this situation
$\phi\in\Irr$ is $2$-large if $\om(\phi)\neq 2^j\om_i$. For
$2$-large representations the estimates from
\cite[Theorem~1.1]{Supl} are available. For remaining
representations certain estimates could be obtained as well, but
this article does not seem a proper place for this. We plan to
handle this question in a subsequent paper which will be devoted to refining
some estimates in \cite{Supl}.

The results of this article as well as those of \cite{Supl} can be easily
transferred to irreducible $K$-representations of finite classical groups
in characteristic $p$.
\section{Notation and preliminary comments}
Throughout the article for a semisimple algebraic group $S$ the
symbols $\Irr(S)$, $\mathbf{X}(S)$, and $\mathbf{X^+}(S)$ mean the
same as the similar ones for $G$ introduced earlier; $R(S)$ is the
set of roots of $S$, $\langle S_1, \ldots, S_j\rangle$ is the
subgroup in $S$ generated by subgroups $S_1, \ldots, S_j$;
$\Irr_p(S)\subset\Irr(S)$ is the set of $p$-restricted
representations, i.e. irreducible representations with
$p$-restricted highest weights; $\mathbf{X}(\phi)$
($\mathbf{X}(M)$) is the set of weights of a representation
$\varphi$ (a module $M$); $\dim M$ is the dimension of $M$;
$M(\omega)$ is the irreducible $S$-module with highest weight
$\omega$; $L$ is the Lie algebra of $G$; $R=R(G)$, $R^+\subset R$
is the set of positive roots; $\Irr_p=\Irr_p(G)$;
$\mathcal{X}_{\beta}\subset G$ and $X_{\beta}\in L$ are the root
subgroup and the root element associated with $\beta\in R$,
$\mathcal{X}_{\pm i}=\mathcal{X}_{\pm\alpha_i}$, and  $X_{\pm
i}=X_{\pm\alpha_i}$. Set
$H(\beta_1,\ldots,\beta_j)=\langle\mathcal{X}_{\beta_1},
\mathcal{X}_{-\beta_1} \ldots, \mathcal{X}_{\beta_j},
\mathcal{X}_{-\beta_j}\rangle$. For  $\omega\in\mathbf{X}(S)$ and
$\alpha\in R(S)$  denote by $\langle\omega,\alpha\rangle$ the
value of the weight $\omega$ on the root $\alpha$. For an
$S$-module $M$ and a unipotent element $x\in S$ define $k_M(x)$
similarly to $k_{\phi}(x)$. If $\abs{x}=p$, then $n_{\phi}(x)$ is
the number of Jordan blocks of size $p$ of the matrix $\phi(x)$
for a representation $\phi$ of $S$ and $n_M(x)$ denotes the same
number for a module $M$ affording  $\phi$.

 An element $x\in G$ of order $p$  can be
 embedded into a closed connected subgroup $\Gamma$ of type $A_1$ whose
  labelled diagram coincides with $\Delta_x$ (see \cite[Theorem 4.2]{LaTe}).
Set $\mathbf{X}_1=\mathbf{X}(A_1(K))$ (the simply connected group
of this type) and identify $\mathbf{X}_1$ with $\mathbb{Z}$
mapping $a\om_1\in\mathbf{X}_1$ into $a\in\mathbb{Z}$. Then
$\mathbf{X}(\Gamma)$ can be identified with a subset of
$\mathbb{Z}$. The canonical homomorphism $\tau_x$ can be obtained
as the restriction of weights from a maximal torus $T\subset G$ to
a maximal torus $T_1\subset \Gamma$ such that $T_1\subset T$. From
now on we fix the tori $T$ and $T_1$, and all weights and roots of
$G$ and $\Gamma$ are considered with respect to $T$ and $T_1$.
Throughout the text  $\ep_i$ with $1\le i\le r+1$ for $G=A_r(K)$
and $1\le i\le r$ otherwise are weights of the standard
realization of $G$ labelled as in \cite[ch. VIII, \S 13]{B2}. Set
$e_i=\tau_x(\varepsilon_i)$. One can choose $\Gamma$, $T$ and
$T_1$ such that the restriction to $\Gamma$ of the natural
representation of $G$ is a direct sum of irreducible components
with $p$-restricted highest weights (see comments in \cite[Section
3]{Te});
 $e_i\geq e_j$ for $i<j$; $e_i\geq0$ if $G=A_r(K)$ and
$i\leq(r+1)/2$; and $e_i\geq0$ for all $i\le r$ if $G\neq A_r(K)$.
If $H\subset G$ is a semisimple subgroup generated by some root
subgroups, then $T_H=T\cap H$ is a maximal torus in $H$. If
$T_1\subset T_H$, we denote by the same symbol $\tau_x$ the
homomorphism $\mathbf{X}(H)\rightarrow\mathbb{Z}$ determined by
restricting weights from $T_H$ to $T_1$. This causes no confusion.
If an element  $v$ of  some $G$-module is an eigenvector for $T$,
we denote its weights with respect to $T$, $T_H$, and $T_1$ by
$\omega(v)$, $\omega_H(v)$, and $\omega_{\Gamma}(v)$. In what
follows $x$ is conjugate to an element of $G_m$, $\abs{x}=p$, $m$
and $r-m$ are such as in  the assertion of Theorem~\ref{mt}, and
$\delta_i=\tau_x(\alpha_i)$, $1\le i\le r$.
\begin{lemma}\label{lc} Set $l=[(m+2)/2]$ for $G=A_r(K)$ and $l=m$
otherwise. Then $c_x=\sum^l_{i=1}\delta_i$ and
 $c_x\leq p-1$.
\end{lemma}

\begin{proof} Put $k=l-1$ for $G=A_r(K)$, $m=2t$, and $k=l$ for $G=A_r(K)$,
$m=2t+1$.  Our assumptions on $e_i$, $m-r$, and $x$ imply that
$e_i=0$ for  $k<i<r+2-k$ if $G=A_r(K)$ and $e_i=0$ for $i>m$
otherwise; notice that $e_{k+1}=e_{k+2}=0$ for $G=A_r(K)$. Now it
follows from the definition of $c_x$ and the formulae in \cite[ch.
VIII, \S 13]{B2}  that
$c_x=\sum^l_{i=1}\delta_i=e_1-e_{l+1}=e_1$. As $e_1$ is a
weight of a $p$-restricted $\Gamma$-module, we have $e_1<p$. This
yields the lemma.
\end{proof}

\begin{proof}[Proof of Theorem~\ref{mt}] Set
$\omega=\omega(\varphi)$ and let $\omega=\sum^r_{i=1}a_i\omega_i$.
It is clear that $\om\neq0$ as $\tau_x(\om)\neq0$. Define
subgroups $H_1$ and $H_2\subset G$ as follows. For $G=A_r(K)$ set
$u=r-t+2$ if $m=2t$ and $r-t+1$ if $m=2t+1$,
$\beta=\varepsilon_{t+1}-\varepsilon_u$,
  \begin{equation*}
H_1=H(\alpha_1,\ldots,\alpha_t, \beta, \alpha_u, \ldots,
\alpha_r), \quad H_2=H(\alpha_{t+2}, \ldots, \alpha_{u-2})
 \end{equation*}
(we have $H_1=H(\alpha_1, \varepsilon_2-\varepsilon_{r+1})$ for
$m=2$ and $H_1=H(\varepsilon_1-\varepsilon_{r+1})$ for $m=1$). For
$G=B_r(K)$, $C_r(K)$, or $D_r(K)$ put $\beta=\varepsilon_m$,
$2\varepsilon_m$, or $\varepsilon_{m-1}+\varepsilon_m$,
respectively, and
 \begin{equation*}
 H_1=H(\alpha_1,\ldots,\alpha_{m-1}, \beta).
 \end{equation*}
Next, set
\[H_2=H(\alpha_{m+1},\dotsc,\alpha_{r-1},\varepsilon_{r-1}+\varepsilon_r)\]
for $G=B_r(K)$ and
\[H_2=H(\alpha_{m+1},\dotsc,\alpha_r)\]
for $G=C_r(K)$ or $D_r(K)$
(here $H_1=H(\beta)$ for $G=C_r(K)$ and $m=1$).
 One easily observes that the sets of roots in brackets used to
 define $H_1$ and $H_2$ yield bases of the systems $R(H_1)$ and
 $R(H_2)$, respectively. Denote these bases by $\mathcal{B}_i$.  In all
cases $H_1$ is conjugate to $G_m$ in $G$.  We have
  $H_2\cong A_{r-m-1}(K)$, $D_{r-m}(K)$,
$C_{r-m}(K)$, or $D_{r-m}(K)$ for $G=A_r(K)$, $B_r(K)$, $C_r(K)$,
or $D_r(K)$, respectively.  It is clear
  that  the
 subgroups $H_1$ and $H_2$ commute. Set $H=H_1H_2$. Let
 $U_i=\langle\mathcal{X}_{\gamma}\mid \gamma\in R^+, \quad
 \mathcal{X}_{\gamma}\subset H_i\rangle$, $i=1,2$, and $U=U_1U_2$. It is
 not difficult to conclude that $U_i$ is a maximal unipotent
 subgroup in $H_i$ and  $U$ is such a subgroup in $H$. We can assume
 that $x\in U_1$, $\Gamma\subset H_1$ and $T_1\subset T_{H_1}$.
 We shall write a weight $\mu\in\mathbf{X}(H)$  in the form
 $(\mu_1, \mu_2)$ where $\mu_i\in\mathbf{X}(H_i)$ is the
 restriction of $\mu$ to $T_{H_i}$.  Set $M=M(\om)$.

It is clear that $n_V(x)=\dim(x-1)^{p-1}V$ for each $H$-module
$V$. Taking this into account, it is not difficult to conclude the
following. If $0\subset W_1\subset\ldots\subset W_t=V$ is a
filtration of $V$, $F_i=W_i/W_{i-1}$, $1\le i\le t$, and
$n_{F_i}(x)=n_i$, then
\begin{equation}\label{eqx}
n_V(x)\geq \sum^t_{i=1}n_i.
\end{equation}

First suppose that  $\phi\in\Irr_p$. Since passing to the dual
representation does not influence the Jordan form of $\varphi(x)$,
one can assume that $a_i\neq0$ for some $i\leq (r+1)/2$ if
$G=A_r(K)$. As for $p$-large representations the estimates of
\cite[Theorem~1.1]{Supl} hold; we also assume that $\phi$ is not
$p$-large. Hence $\langle\mu, \alpha\rangle<p$ for all
$\mu\in\mathbf{X}(\phi)$ and long roots $\alpha$ (for all $\alpha$
if $G=A_r(K)$ or $D_r(K)$). By the formulae for the maximal roots
of the classical groups in \cite[Tables 1-4]{B1}, this forces that
\begin{equation}\label{eqp}
\begin{split}
a_1+\ldots+a_r<p \quad  & \text{for}\quad
G=A_r(K)\quad\text{or}\quad C_r(K), \\
 a_1+2a_2+\ldots+2a_{r-1}+a_r<p \quad & \text{for} \quad
 G=B_r(K),\\
 a_1+2a_2+\ldots+2a_{r-2}+a_{r-1}+a_r<p \quad & \text{for} \quad
G=D_r(K).
\end{split}
\end{equation}
Now we proceed to construct
 two composition factors $M_1$ and $M_2$ of the restriction
 $M|H$ such that $n_{M_1}(x)\geq d(r-m)$ and $n_{M_2}(x)>0$.
This will be done for almost all $\omega$.
 In exceptional cases we shall find  one  factor $M_1$ such that
 $n_{M_1}(x)>d(r-m)$. By (\ref{eqx}), this would yield
 the assertion of the theorem.

 Let $v\in M$ be a nonzero highest weight
 vector. Put $\mu_i=\omega_{H_i}(v)$. The vector $v$ generates an
 indecomposable $H$-module $V_1$ with highest weight
 $\mu=(\mu_1, \mu_2)$. Using (\ref{eqp}), one can deduce that
 $\langle\mu_1, \beta\rangle<p$ for all $\beta\in\mathcal{B}_1$.
 Here for $G=B_r(K)$ we take into account that $m>1$. Hence
 $\mu_1$ is $p$-restricted. Now assume that either $G\neq B_r(K)$, or
 $a_i\neq0$ for some $i<r$. For such representations  we construct another
  weight vector $w\in M$ that is fixed by $U$. Set $l=t+1$ for $G=A_r(K)$,
$m=2t$; otherwise take $l$ as in
Lemma~\ref{lc}. First suppose that $a_j\neq0$ for some $j\le l$ (Case 1).
Choose maximal such $j$ and put $w=X_{-l}\ldots X_{-(j+1)}X_{-j}v$. Now let
$a_j=0$ for all $j\le l$ (Case 2). Our assumptions on $a_i$ imply that
$a_i\neq0$ for some $i>l$; furthermore, one
can take $i\le(r+1)/2$  for $G=A_r(K)$ and $i<r$ for $G=B_r(K)$. Choose
minimal such $i$ and set $w=X_{-l}\ldots X_{-(i-1)}X_{-i}v$ if $G\neq D_r(K)$
 or $i<r$ and $w=X_{-l}\ldots X_{-(r-3)}X_{-(r-2)}X_{-r}v$ for $G=D_r(K)$ and
$i=r$. It follows from \cite[Lemma 2.1(iii) and Lemma 2.9]{Supl} that in all
cases $w\neq0$. Using
\cite[Lemma~72]{St-b} and analyzing the roots in $\mathcal{B}_1$
and $\mathcal{B}_2$ and the weight system $\mathbf{X}(\phi)$, we
get that $U$ fixes $w$ in all situations. Here it is essential
that the case $G=B_r(K)$ with $\om=a_r\om_r$ is excluded. In the
latter case we cannot assert that $\mathcal{X}_{\beta}$ fixes $w$.
Set $\lambda_i=\om_{H_i}(w)$, $i=1,2$. Now it is clear that $w$
generates an indecomposable $H$-module $V_2$ with highest weight
$\lambda=(\lambda_1, \lambda_2)$. We claim that $\lambda_1$ is
$p$-restricted. Write down all the situations where
$\langle\lambda_1, \gamma\rangle\neq \langle\mu_1, \gamma\rangle$
for some $\gamma\in\mathcal{B}_1$.  We have
$\langle\lambda_1, \beta\rangle=\langle\mu_1, \beta\rangle-1$ in Case 1 if
$j=l$ and $G\neq B_r(K)$ or $j=l-1$ and $G=D_r(K)$ and in Case 2 for
$G\neq B_r(K)$ and all $i$; and $\langle\lambda_1, \beta\rangle=\langle\mu_1,
\beta\rangle-2$ for $G=B_r(K)$ both in Case 1 with $j=l$ and in Case 2. In
Case 1 we also have
$\langle\lambda_1,\alpha_{j-1}\rangle=\langle\mu_1, \alpha_{j-1}\rangle+1$
if $j>1$  and
$\langle\lambda_1, \alpha_j\rangle=\langle\mu_1,\alpha_j\rangle-1$
if $j<l$. In Case 2 one gets
$\langle\lambda_1,\alpha_{l-1}\rangle=\langle\mu_1, \alpha_{l-1}\rangle+1$
if $l>1$.
In all other situations we have
$\langle\lambda_1, \gamma\rangle=\langle\mu_1, \gamma\rangle$. Now
apply (\ref{eqp}) to conclude that $\lambda_1$ is $p$-restricted.

Set $M_1=M(\mu)$, $M_2=M(\lambda)$, $M^j_1=M(\mu_j)$, and
$M^j_2=M(\lambda_j)$, $j=1,2$. Obviously, $M_i$ is a composition
factor of $V_i$. It is well known that $M_i=M^1_i\otimes M^2_i$.
It is clear that $\tau_x(\mu_1)=\tau_x(\omega)\geq p$. Since $x\in
H_1$, we have $\delta_t=0$ if $\alpha_t\in\mathcal{B}_2$.
So by Lemma~\ref{lc},
 \begin{equation*}
\tau_x(\lambda_1)=\tau_x(\omega(w))\geq
\tau_x(\omega)-\sum^l_{i=1}\delta_i=\tau_x(\omega)-c_x\geq p-1.
 \end{equation*}
 It follows from \cite[Theorem~1.1, Lemma~2.5,  and Proposition~2.12]{Sump}
that $k_{M^1_i}(x)=p$. Hence $n_{M_i}(x)\geq  \dim M^2_i$.  One
easily observes that $M^2_1$ and $M^2_2$ cannot both be trivial
$H_2$-modules. Our assumptions on $r-m$ and
\cite[Proposition~5.4.13]{KL} imply that the dimension of a
nontrivial irreducible $H_2$-module is at least $d(r-m)$. In the
exceptional case where $G=B_r(K)$ and $\om=a_r\om_r$ we need to
evaluate $\dim M^2_1$.
First let $a_r\not=1$. As above, $X_rv\not=0$. This implies that
$\mathbf{X}(M_1^2)$ contains a dominant weight $\mu_2-\alpha_r$ and $\dim
M_1^2$ is greater than the size of the orbit of $\mu_2$ under the action of the
Weyl group of $H_2$. The latter is equal to $2^{r-m-1}\ge d(r-m)$ for our
values of $r-m$. By (1), this yields the assertion of the theorem for almost
all $\varphi\in\operatorname{Irr}_p$. It remains to consider the case where
$G=B_r(K)$ and $\omega=\omega_r$. It is well known that then the restriction
$M\downarrow H_1$ is a direct sum of $2^{r-m}$ $H_1$-modules $N=M(\omega_m)$.
Since $k_m(x)=p$, we get $k_N(x)=p$ and $n_M(x)\ge 2^{r-m}>d(r-m)$.

Now suppose that $\phi\in\Irr\backslash\Irr_p$. By the Steinberg
tensor product theorem \cite[Theorem~1.1]{St}, $\phi$ can be
represented in the form $\bigotimes^s_{j=1}\phi_j\Fr^j$ where $\Fr$
is the Frobenius morphism of $G$ associated with raising elements
of $K$ to the $p$th power and all $\phi_j\in\Irr_p$. It is clear
that the morphism $\Fr$ does not influence the Jordan form of
$\phi(x)$. Hence one can assume that $\phi=\psi\otimes\theta$
where $\theta=\phi_j\Fr^j$ for some $j$ and both $\psi$ and
$\theta$ are nontrivial. Set $a=\sigma_x(\om(\psi))$,
$\nu=\om(\phi_j)$, $b=\tau_x(\nu)$ and define by $\mu$ the
restriction of $\nu$ to $T_H$. Now it follows from the definitions
of $\sigma_x$ and $\tau_x$ that $\sigma_x(\om)=a+b$. By
\cite[Theorem~1.1, Lemma~2.5 and Proposition~2.12]{Sump},
$k_{\psi}(x)=\min\{a+1,p\}$ and $k_{\theta}(x)=\min\{b+1,p\}$.
First suppose that $a$ or $b\geq p-1$. Set $\rho=\psi$ if $a\geq
p-1$ and $\rho=\theta$ otherwise and denote by $\pi$ the remaining
representation from the pair $(\psi, \theta)$. Then
$k_{\rho}(x)=p$ and  \cite[ch. VIII, Lemma~2.2]{Fe}  implies that
$n_{\phi}(x)\geq \dim\pi$. Let $d(r)$ be the value of $d(r-m)$ if
one formally sets $m=0$. Then by \cite[Proposition~5.4.13]{KL},
$\dim\pi\geq d(r)>d(r-m)$ which settles the case under
consideration.

Now assume that both $a$ and $b<p-1$. Then $k_{\psi}(x)=a+1$ and
$k_{\theta}(x)=b+1$. Since $\sigma_x(\om)\geq p-1+c_x$, we have
$b>c_x$. Arguing as for $p$-restricted $\phi$, we can and shall suppose
that $\langle\nu, \alpha_i\rangle\neq0$ for some $i\leq (r+1)/2$ if
$G=A_r(K)$. Put $M'=M(\nu)$ and construct the  composition factors
$M_i$, $i=1,2$, of the restriction $M'|H$ as for  $p$-restricted
$M$ before. Transfer the notation $\mu_1$, $\lambda_1$, and
$M^i_j$, $i,j=1,2$, to $M'$. Again we have the exceptional case
$G=B_r(K)$ and $\nu=a_r\om_r$ where we do not construct $M_2$ and
consider $M_1$ only. Obviously, $\tau_x(\mu_1)=b$. As before, we
deduce that $\tau_x(\lambda_1)\geq b-c_x$. By
 \cite[Theorem~1.1, Lemma~2.5, and Proposition~2.12]{Sump}, $k_{M_1}(x)=b+1$ and
$k_{M_2}(x)\geq b+1-c_x$.
Let $n_i$ be the number of Jordan blocks of the maximal size in the canonical
form of $x$ as an element of $\operatorname{End} M_i$, $i=1,2$. Looking at the
realizations of $M_i$ as tensor products, one easily observes that $n_i\ge\dim
M_i^2$.
Set $F_1=M(\omega(\psi))\otimes M_1$,
$F_2=M(\omega(\psi))\otimes M_2$ and consider $F_i$ as $H$-modules
in the natural way. In the general case the $H$-module $M$ has a
filtration two of whose quotients are isomorphic to $F_1$ and
$F_2$, respectively. In the exceptional case $F_1$ is a quotient
of a submodule in $M$. Observe that $a+k_{M_2}(x)\geq p$. Using
\cite[ch. VIII, Theorem~2.7]{Fe} that describes the canonical
Jordan form  of a tensor product of unipotent blocks, we obtain
that $k_{M_i}(x)=p$ and $n_{F_i}(x)\geq \dim M^2_i$. As for
$p$-restricted $M$, we show that $n_1\ge 2^{r-m}$ if $G=B_r(K)$ and
$\nu=\omega_r$ and conclude that $\dim M_1^2+\dim M_2^2>d(r-m)$ in the general
case and $\dim\,M_1^2>d(r-m)$ in the exceptional cases with $a_r\not=1$.
Now
(\ref{eqx}) completes the proof. \end{proof}
\begin{proof}[Proof of Proposition~\ref{mp}] Let $a$, $x$, $m$, and
$c$ be such as in the assertion of the proposition. Assume that
$p<ac\leq p+c-1$. Therefore
 we have $(a-1)c\leq p-1$. Set $M=M(\omega)$ and denote by
$M_t$ the weight subspace of
weight $t\in\mathbb{Z}$ in the $\Gamma$-module $M$. It is clear
that the Weyl group of $\Gamma$ interchanges $M_t$ and $M_{-t}$;
hence $\dim M_t=\dim M_{-t}$. Put  $e=(a-1)c$,
$V_1=\bigoplus_{t>e}M_t$, $V_2=\bigoplus_{t>e}M_{-t}$, and $V=M_e$. Set
$f=[(m+1)/2]$ for $G=A_r(K)$, $f=m$ for $G=B_r(K)$ or $C_r(K)$,
and $f=m-1$ for $G=D_r(K)$. Let $v\in M$ be a nonzero highest
weight vector and put $w=X_{-f}\ldots X_{-2}X_{-1}v$. By
\cite[Lemma 2.9]{Supl}, $w\neq0$.  We need a
subgroup $S$ which can be defined as follows. Put $I=\{i\mid
1\le i\le r,\quad \delta_i=0\}$ and $S=\langle\mathcal{X}_i,
\mathcal{X}_{-i}\mid i\in I\rangle$. The canonical Jordan forms of $x$
in the standard realizations of $G_m$ and $G$ are well known. We have
$\abs{x}=p$ since the dimension of the first realization is at most $p$
due to our assumptions. Taking into
account these Jordan forms, one easily obtains the values of $\delta_i$,
$1\le i\le r$, and using Lemma~\ref{lc}, deduces the following facts:
$I=\{i\mid f+1\leq i\leq r-f\}$ for $G=A_r(K)$ and
$m=2f$, $I=\{i\mid f+1\leq i\leq r\}$ for $G=B_r(K)$ and $D_r(K)$,  and $S=H_2$
in all other cases where $H_2$ is the subgroup defined in the
proof of Theorem~\ref{mt}; $c_x=\sum^f_{i=1}\delta_i=c$,
$\tau_x(\omega)=ac$; and $w\in V$.   Next, observe that  $S\cong A_{r-m}$ for
$G=A_r(K)$ and $m=2f$, $S\cong B_{r-m}(K)$ for $G=B_r(K)$,
and $S\cong D_{r-m+1}$ for $G=D_r(K)$. Our
construction of the vector $w$ shows that $\mathcal{X}_i$ fixes
$w$ if $i\in I$. This forces that $w$ generates an indecomposable
$S$-module $M_S$ with highest weight $\omega_S(w)$. Then one
immediately concludes that $M_S\cong M(\omega_1)$. This
 yields that $\dim M_S=r-m+1=d(r-m)+1$ for $G=A_r(K)$ and $m=2f$,
 $\dim M_S=2(r-m)+1=d(r-m)+1$ for $G=B_r(K)$,
$\dim M_S=2(r-m+1)=d(r-m)+2$ for $G=D_r(K)$, and $\dim M_S=d(r-m)$
 otherwise. It is clear that $M_S\subset V$. Denote by
$\mathbf{X}_f\subset\mathbf{X}(M)$ the subset of weights of the
form $\omega-\sum^f_{i=1}b_i\alpha_i$ and by $M_A$ the irreducible
$A_f(K)$-module with highest weight $a\omega_1$. By Smith's theorem
\cite{Sm}, for each $\mu\in\mathbf{X}_f$ the dimension of the
weight subspace $M_{\mu}\subset M$ coincides with that of the
weight subspace in $M_A$ whose weight differs from $a\omega_1$ by
the same linear combination of the simple roots. Hence $\dim
M_{\mu}$ does not depend upon $r$. Set
$W=\bigoplus_{\mu\in\mathbf{X}_f}M_{\mu}$. Since $M$ is an
irreducible $L$-module and $p>2$, observe that $M$ is a linear
span of vectors of the form $X_{-i_s}\ldots X_{-i_2}X_{-i_1}v$.
Now, analyzing the weight structure of $M$, we conclude that
$V_1\subset W$ and $V=(V\cap W)\oplus M_S$. This implies that
$\dim V_1$ ($=\dim V_2$) and $\dim (V\cap W)$ do not depend upon
$r$.

 It follows from \cite[Lemma 72]{St-b} that
 \begin{equation}\label{eqpo}
 (x-1)^{p-1}M_t\subset\bigoplus_{i\geq t+2p-2} M_i.
 \end{equation}
  Let $M_t\not\subset V_2$. Then $t\geq -e$. Obviously, $e<p-1$ if
  $ac<p-1+c_x$ and $e=p-1$ for $ac=p-1+c_x$. Thus (\ref{eqpo}) implies
  that
  \begin{equation*}
  (x-1)^{p-1}M_t\subset\bigoplus_{t>p-1}M_t\subset V_1
  \end{equation*}
  in the first case and
  \begin{equation*}
(x-1)^{p-1}M_t\subset\bigoplus_{t\geq p-1}M_t\subset V_1\oplus V
\end{equation*} in the second case. This forces that
 $n_M(x)\leq\dim V_2+\dim V_1=2\dim V_1$ in the first case and
$n_M(x)\leq 2\dim V_1+\dim(V\cap W)+\dim M_S$ in the second case.
We have seen before that $\dim M_S=d(r-m)+u$ with $u=0$, $1$, or
$2$. Hence one can take $N_G(a,m,p)=2\dim V_1$ and
$Q_G(a,m,p)=2\dim V_1+\dim(V\cap W)+u$ to complete the proof.
\end{proof}

\begin{remark} For $G=A_r(K)$ or $C_r(K)$
we could give a shorter proof of Proposition~\ref{mp} using the
realization of $\phi$ in the $a$th symmetric power of the standard
module (see \cite[1.14 and  8.13]{Se}), but we need the proof
above for $B_r(K)$ and $D_r(K)$. \end{remark}
 \bigskip

This research  has been supported by the Institute of Mathematics
of the National Academy of Sciences of Belarus in the framework of
the State program ``Mathematical structures'' and by the Belarus
Basic Research  Foundation, Project F\,98-180.

\end{document}